\documentclass[reqno, 11pt]{amsart}

\usepackage{hyperref}
\usepackage{amsmath}
\usepackage{bbm}
\usepackage{amssymb}
\usepackage{amsthm}
\usepackage{xcolor}
\usepackage{orcidlink}

\newtheorem{thm}{Theorem}

\newtheorem{prop}{Proposition}
\newtheorem{lem}{Lemma}

\title{Generalized Foguel-Hankel operators}

\author{Nikolaos Chalmoukis \orcidlink{0000-0001-5210-8206}
} \address{Dipartimento di Matematica e Applicazioni, Universit\'a degli studi di Milano Bicocca, via Roberto Cozzi, 55 20125, Milano, Italy}
\email{nikolaos.chalmoukis@unimib.it}
\thanks{Both authors are members of the INdAM group GNAMPA and are partially supported by the grant INdAM-GNAMPA Project, CUP E53C25002010001, "Transferring Harmonic Analysis between Discrete Structures and Manifolds". The first author was partially supported by the research grant "Yields of the ubiquity and the geometry of inner
functions (YoungInFun)", PID2024-160326NA-I00.}

\author{Giovanni Marano \orcidlink{0009-0001-9082-7882}}

\address{Dipartimento di Matematica e Applicazioni, Universit\'a degli studi di Milano Bicocca, via Roberto Cozzi, 55 20125, Milano, Italy}
\email{g.marano7@campus.unimib.it}
\date{}

\subjclass{47B02, 47B35, 47B32}
\keywords{Foguel-Hankel operators, Kreiss condition, Power boundedness, Polynomial boundedness, Similarity to a contraction}

\begin{document}

\begin{abstract}
	In this paper we introduce a more general class of Foguel-Hankel operators, where the unilateral shift on $\ell^2(\mathbb{N}) $  is replaced by a general multiplication operator on the Hardy space $H^2$ . We prove that Peller's condition is sufficient for the operator to be power bounded, but in general it is not necessary. When the Hankel matrix is the Hilbert matrix, we prove that being similar to a contraction is equivalent to the (a priori) weaker Kreiss condition. 
\end{abstract}

\maketitle

\section{Introduction}

The starting point of this work is the r\^ole played by Foguel-Hankel operators 
as a rich source of examples in the family of problems related to similarity to a contraction. 
A bounded linear operator $T $ acting on a Hilbert space $H $ is called {\it similar to a contraction} if there exists another bounded invertible operator $ M $ and a contraction $K $ such that 
\begin{equation}\tag{SimC} T = M^{-1}KM. \end{equation}

It follows directly from von Neumann's inequality \cite{Neumann} that if $T $ is similar to a contraction it is also {\it polynomially bounded}, that is, there exists $C>0 $ such that 
\begin{equation}\tag{PolB}
	\Vert p(T)\Vert \leq C \sup_{z\in \mathbb{D}}|p(z)|,
\end{equation}
for every polynomial $p $, where $\mathbb{D} $ is the unit disc. 
In particular, by considering only the monomials, a polynomially bounded operator $T $ has powers which are uniformly bounded in norm, i.e., there exists $C>0 $ such that 
\begin{equation}
	\tag{PB}
	\Vert T^n\Vert\leq C, \,\,\, \forall n \in \mathbb{N}=\{1,2,3,\dots \}.
\end{equation}

An even weaker condition is the so-called Kreiss condition which requires that the spectrum of the operator $T $ is contained in the closed unit disc and there exists $ C>0$ such that 
\begin{equation}
	\tag{Kreiss}
	(|\lambda|-1)\Vert (\lambda-T)^{-1}\Vert \leq C, \,\,\, \forall |\lambda|>1.
\end{equation}

It turns out that none of the above implications is reversible, but this is far from obvious.
An operator that satisfies the Kreiss condition but is not power bounded has been known for a long time \cite{McCarthy1965} (see also \cite{Chalmoukis2025} for a simple example).
The counterexamples to the other two implications are operators of Foguel-Hankel type.

Let us now define what a (scalar) Foguel-Hankel operator is. 
Consider the Hardy spaces of the unit disc $H^p(\mathbb{D}), 1\leq p<\infty $  consisting of holomorphic functions in $\mathbb{D} $  such that 
\[ \Vert  f\Vert_{p}^p := \sup_{0<r<1}\frac{1}{2\pi} \int_0^{2\pi}|f(re^{i\theta})|^p d\theta < + \infty. \]
 
It is well known \cite{Garnett2006} that we can identify isometrically $H^p(\mathbb{D}) $ with the subspace $H^p(\mathbb{T}) $  of $L^p( \mathbb{T} ) $, where $\mathbb{T} = \partial \mathbb{D} $, which consists of functions having vanishing negative Fourier coefficients. We adopt this convention throughout and we will write simply $H^p $ to denote any of the two spaces when confusion does not arise.   
We now consider some operators which act on the Hilbert space $H^2(\mathbb{T}) $. 
Let $J $ be the {\it flip} operator which is defined by the rule $Jf(z)=f(\overline{z}) $ for $f\in L^2(\mathbb{T}) $. 
We will denote by $P_- $ the orthogonal projection from $L^2(\mathbb{T}) $ onto the space $H^2_-(\mathbb{T}):= J(H^2(\mathbb{T})) $. 
Given now $f\in H^2(\mathbb{T}) $ we can define the {\it Hankel operator} $H_f $  with symbol $f $ on the dense subspace of analytic polynomials as follows  
\[ H_f(p) := P_-(J(f)p).  \]
 Nehari's theorem gives a description of the functions $f $ for which $H_f$ extends to a bounded operator on the whole $H^2(\mathbb{T}) $. In fact, we have that 
 \[ \Vert H_f \Vert = \sup_{u\in H^2, \,\, \Vert u\Vert_{1}\leq 1} |\langle f, u \rangle|, \]
where the scalar product is the one associated with the Hardy space norm. The space of all holomorphic functions $f\in H^2 $ for which the above quantity is finite will be denoted by $BMOA $.   
 Another well studied operator acting on the Hardy space is the unilateral shift $S $, which acts by multiplication by the independent variable $z $. 
 We now have the ingredients needed in order to define a Foguel-Hankel operator $\Gamma_f $ which acts on the direct sum  $ H^2(\mathbb{T}) \oplus H^2(\mathbb{T})$  of two Hardy spaces equipped with the standard Hilbert space norm.  $\Gamma_f $ is defined as the following $2\times 2 $ block matrix
 \[ \Gamma_f = \begin{bmatrix} S^* & JH_f \\ 0 & S \end{bmatrix}. \]
 
 A key feature of this class is that its polynomial functional calculus can be expressed in a very simple way. Let $p\in \mathbb{C}[z] $ be a polynomial, then as  can be readily verified 
 \begin{equation} \label{eq:FHpol} p(\Gamma_f) = \begin{bmatrix} p(S^*) & JH_{p'(S^*)f} \\ 0 & p(S) \end{bmatrix}. \end{equation}
 Thus polynomial expressions in $\Gamma_f $ reduce to a block matrix containing multiplication operators by $p $ and its adjoint and another Hankel operator, which makes estimates particularly accessible. 

 With respect to the spectral conditions discussed above, in \cite{Peller1982}, Peller showed that the operator $\Gamma_f $ is power bounded if and only if $f' $ is a Bloch function, that is 
 \begin{equation*}
 \Vert f'\Vert_\mathcal{B}:=	\sup_{z\in \mathbb{D}}|f''(z)|(1-|z|^2) < +\infty. 
 \end{equation*}
 
 We will denote by $\mathcal{B} $ the space of Bloch functions. In the papers \cite{Peller1983, Aleksandrov1996} the authors showed that $\Gamma_f $ is polynomially bounded if and only if $f'$ belongs to the space $  BMOA $. Consequently, there exists a Foguel-Hankel operator which is power bounded but not polynomially bounded. 

 Regarding similarity to a contraction, it turns out (see \cite{Bourgain1986}) that for the Foguel-Hankel operators considered here, similarity to a contraction is equivalent to polynomial boundedness. An example of a polynomially bounded operator which is not similar to a contraction was constructed in a celebrated work by Pisier \cite{Pisier1996}.

Our first result gives the missing implication needed to compare the Kreiss condition with Peller's power boundedness criterion for classical Foguel-Hankel operators.

 \begin{thm}\label{thm:FH_Kreiss}
 Let $\Gamma_f $ be a Foguel-Hankel operator which satisfies the Kreiss condition. Then $f' $ belongs to the Bloch space. 	
 \end{thm}

 Combining Theorem \ref{thm:FH_Kreiss} with Peller's characterization of power boundedness, we obtain that, in the classical Foguel-Hankel class, the Kreiss condition is equivalent to power boundedness. 

 Motivated by this rich source of examples, in this work we attempt to define a more general class of operators that we think are worth investigating. Consider a bounded analytic function $\varphi $ and let $M_\varphi: H^2(\mathbb{D}) \to H^2(\mathbb{D}) $ be the operator of multiplication by $\varphi $, $M_\varphi f = \varphi f. $ 
 It is a standard fact that $ \Vert M_\varphi \Vert = \Vert \varphi\Vert_\infty = \sup_{|z|<1}|\varphi(z)|. $ 
 We will use the notation $ \widetilde{\varphi}(z) = \overline{\varphi(\overline{z})}$. Let also $f\in BMOA $ so that $H_f $ is bounded on $H^2 $. The generalized Foguel-Hankel operator $\Gamma_{f,\varphi} $  with symbols $f $ and $\varphi $ is defined as the $2\times 2 $ block operator acting on $H^2 \oplus H^2 $ given by  
\begin{equation} \label{def:GFH} \Gamma_{f,\varphi} = \begin{bmatrix} M_{\widetilde{\varphi}}^* & JH_f \\ 0 & M_\varphi  \end{bmatrix} : H^2 \oplus H^2 \to H^2 \oplus H^2. \end{equation}
\noindent Observe that the upper triangular form of the operator $\Gamma_{f,\varphi} $ implies that 
\begin{equation*}
	\sigma(M_\varphi) \subseteq \sigma(\Gamma_{f,\varphi}) \subseteq \sigma(M_{\widetilde{\varphi}}^*) \cup \sigma(M_\varphi),
\end{equation*}
where $\sigma $ denotes the spectrum of an operator, which implies that 
\[ \sigma(\Gamma_{f,\varphi}) = \overline{\varphi(\mathbb{D})}. \]

Hence, since all the conditions that are relevant for us imply that $ \sigma(\Gamma_{f,\varphi}) \subseteq \overline{\mathbb{D}} $, from now on we will always require that $\varphi $ is a holomorphic self-map of the unit disc.

This family of operators preserves some of the nice algebraic properties of the classical Foguel-Hankel operators $\Gamma_f $. 
A key property is the existence of a similarly simple expression for polynomials of the operator $\Gamma_{f,\varphi} $. 
Let $p\in \mathbb{C}[z] $ be a polynomial, then (see Lemma \ref{lem:FHAlg}) we have that  

\begin{equation} \label{eq:polcalc}
	p(\Gamma_{f,\varphi}) = \begin{bmatrix} M^*_{\widetilde{p\circ \varphi}} & JH_{M^*_{\widetilde{p'\circ \varphi}}f} \\ 0 & M_{p \circ \varphi} \end{bmatrix} = \Gamma_{M^*_{\widetilde{p'\circ \varphi}  }f, p\circ \varphi}.
\end{equation}

\noindent So the class of generalized Foguel-Hankel operators is closed under polynomial expressions. 

Concerning power boundedness for the operator $ \Gamma_{f,\varphi} $ the next theorem shows that, independently of the behaviour of $\varphi $ near the boundary, the generalized Foguel-Hankel operator cannot behave worse than its classical counterpart.  

\begin{thm} \label{thm:FHPb}
	Let $f' $ be a Bloch function and $\varphi$ a holomorphic self-map of the unit disc. Then the operator $\Gamma_{f,\varphi} $ is power bounded.  
\end{thm}

One might hope that the equivalence between the Bloch condition and power boundedness continues to hold for arbitrary  $f\in BMOA $,  but the situation seems to be more complicated. 
We have not been able to give a complete characterization of all pairs  $(f,\varphi)$ such that the operator is power bounded, a concrete counterexample and a more extensive discussion around this problem can be found in Section \ref{sec:fixed_symbol}.
In a different direction, consider the symbol $ f(z)=\sum_{n=0}^{\infty} \frac{z^n}{n+1} $ so that the corresponding Hankel operator is the so-called Hilbert matrix $\mathcal{H} $, which is represented with respect to the standard orthonormal basis of $H^2 $ and $H^2_- $  by the matrix $ [ (1+n+m)^{-1} ]_{n,m =0}^\infty $.   
In this case, we will refer to the generalized Foguel-Hankel operator as a Foguel-Hilbert operator with symbol $\varphi $ and we will denote it by $h_\varphi $. 
The following theorem holds in this case. 
\begin{thm} \label{thm:Hilbert}
	Let $\varphi : \mathbb{D} \to \mathbb{D} $ be a holomorphic function and $h_\varphi $ the corresponding Foguel-Hilbert operator. The following are equivalent:
	\begin{itemize}
		\item[(i)] The operator $h_\varphi $ is similar to a contraction. 
		\item[(ii)] The operator $h_\varphi $ satisfies the Kreiss condition 
		\item[(iii)] It holds that $\limsup_{r\to 1^-}|\varphi(r)| < 1 $. 
	\end{itemize}
\end{thm}

We therefore see that in this case all the conditions discussed above collapse to a single one. Furthermore, an interesting feature of this theorem is that the behaviour of the operator depends only on the modulus of the function $\varphi $ on the segment $(0,1) $. 

The paper is organized as follows. In Section \ref{sec:preliminaries} we recall some basic results and definitions. 
The tools that we use are largely standard, except perhaps for Luecking's embedding theorems for derivatives of Hardy functions. 
In Section \ref{sec:proofs} we give the proofs of Theorems \ref{thm:FH_Kreiss} and \ref{thm:FHPb}, while in Section \ref{sec:Hilbert} we discuss further the case that the Hankel matrix is the Hilbert operator and we prove Theorem \ref{thm:Hilbert}. Finally, in Section \ref{sec:open_problems} we briefly discuss some problems that currently lie outside of our methods and appear to be of interest.

\section{Preliminaries} \label{sec:preliminaries}
In the rest of this work $m $ will denote the Lebesgue area measure normalized so that $m(\mathbb{D}) = 1 $ and $ | \cdot | $ will denote the arc length measure on $\mathbb{T}. $
Let us begin by discussing some material related to spaces of analytic functions.   
We have already introduced the family of Hardy spaces $H^p, 1\leq p \leq \infty $. 
A space which is closely related to the Hardy spaces is the space of functions of bounded mean oscillation $BMOA$. 
There exist a number of equivalent definitions of $BMOA $ but for our purposes it is defined as the dual of $H^1 $ with respect to the usual Cauchy pairing. 
Explicitly, a holomorphic function $f $ belongs to $BMOA $ if $f\in H^2 $ and 
\begin{equation*}
	\Vert f\Vert_* := \sup_{u\in H^2, \,\, \Vert u\Vert_1 \leq 1 } | \langle f, u \rangle| < + \infty.
\end{equation*}

It will be more convenient to use a Fefferman-Stein type expression for the pairing with respect to the scalar product of the Hardy space. 
The next formula is well known, but we give a short proof for completeness.

\begin{lem} \label{lem:PL}
	Let $f, g \in H^2 $. Then it holds that 
	\[ \langle f, g \rangle_{H^2} = f(0)\overline{g(0)}+f'(0)\overline{g'(0)}+2 \int_{\mathbb{D}}f''(z)\overline{g''(z)}\big( \ln\frac{1}{|z|}(1+|z|^2)-(1-|z|^2) \big)d m(z). \]
	\end{lem}

\begin{proof}
	It suffices to prove the identity for $f=g $ such that $f(0)=f'(0)=0 $ and by taking dilations of $f $ we can also assume that it is analytic in an open neighbourhood of $ \overline{\mathbb{D}}.$   

	Observe that the function $H(z)=\frac{1}{2}\big( \ln\frac{1}{|z|}(1+|z|^2)-(1-|z|^2) \big) $ is positive for $|z|<1 $ and it satisfies 
	\begin{equation}
		H(z)=\frac{\partial H}{\partial r}(z) = 0, \,\, |z|=1, \,\,\, \Delta H(z) = \ln \frac{1}{|z|^2}, 0<|z|<1,
	\end{equation}
	where $\frac{\partial }{\partial r} $ is the radial derivative.  
	Then apply first the Fefferman-Stein formula and then Green's theorem as follows 
	\begin{align*}
		\Vert f\Vert^2_2 & = \int_\mathbb{D} |f'(z)|^2 \ln \frac{1}{|z|^2} d m(z) \\
					   & = \lim_{\varepsilon \to 0 }\int_{\varepsilon <|z| <1} |f'(z)|^2 \Delta H(z) d m(z) \\
					   & =4 \int_{\mathbb{D}} |f''(z)|^2 H(z) d m(z) \\
					   & \quad - \lim_{\varepsilon \to 0} \int_{|z| = \varepsilon} |f'(z)|^2 \frac{\partial H}{\partial r }(z) - H(z)\frac{\partial |f'|^2}{\partial r}(z) \frac{|dz| }{\pi} \\
					   & = 4 \int_{\mathbb{D}} |f''(z)|^2 H(z) d m(z).
	\end{align*}
\end{proof}

The next result we are going to discuss is the classical embedding theorem of Carleson and a generalization of it, due to Luecking.
For a point $z\in \mathbb{D} $ let $D(z) $ denote the disc centered at $z $ with radius $\frac{1}{2}(1-|z|) $. 
Furthermore, given an arc $I\subseteq \mathbb{T} $, we denote by $S(I) $ the Carleson box associated with $I, $ that is  

\[ S(I)= \{ re^{i\theta}\in \mathbb{D} : e^{i\theta}\in I, 1-r \leq |I| \}. \]

Given a positive Borel measure $\mu $ in the unit disc, Carleson's embedding theorem claims that the inequality 
\[ \int_{\mathbb{D}}|u(z)|d \mu(z) \leq A \Vert u\Vert_{1}, \]
holds for all $u \in H^1  $ and for some constant $A>0 $ independent of $u $, if and only if there exists another constant $B>0 $ such that for every arc $I \subseteq \mathbb{T} $  

\[ \mu(S(I)) \leq B  |I|.   \]

\noindent The following generalization is due to Luecking \cite[Theorem 1]{Luecking1991}. 
Notice that in his article, Luecking states the theorem for the Hardy space of the half-plane.  The version we report below is obtained by a conformal change of variables. 

\begin{thm} \label{thm:Lue}
	Let $\mu $ be a positive finite Borel measure in $\mathbb{D} $ and $d\in \mathbb{N} $ . 
	Then, there exists a constant $A>0 $ such that 
	\begin{equation}\label{eq:derCar} \int_\mathbb{D}|u^{(d)}(z)| d\mu(z) \leq A \Vert u\Vert_{1}, \,\,\, \forall u \in H^1, \end{equation}
	if and only if there exists a constant $B $ such that 
	\begin{equation}\label{eq:Lue}
		\Big(\frac{1}{|I|}\int_{S(I)} \frac{\mu(D(z))^2}{(1-|z|)^{2d+3}}  d m(z) \Big)^{\frac{1}{2}} \leq B,
	\end{equation}
	for all arcs $I \subseteq \mathbb{T}. $ 
\end{thm}

In the above theorem the constants $A $ and $B $ depend on $d $ and the measure $\mu $. Moreover, if $A, B $ denote the best possible constants in \eqref{eq:derCar} and \eqref{eq:Lue} respectively, there exists a constant $C $ depending only on $d $ such that $A/C \leq B \leq CA  $. The same, of course, holds in Carleson's embedding theorem.

The next proposition is a generalization of Schwarz Lemma for holomorphic self-maps of $\mathbb{D} $. 
\begin{lem}\label{lem:selfmaps}
	Let $\varphi : \mathbb{D} \to \mathbb{D} $ be a holomorphic function. 
	Then the following is true:
	\begin{itemize}
		\item[(i)] $(1-|z|^2) \leq \frac{1+|\varphi(0)|}{1-|\varphi(0)|}(1-|\varphi(z)|^2) $,
		\item[(ii)] $(1-|z|^2)^d|\varphi^{(d)}(z)| \leq d! (1-|\varphi(z)|^2)(1+|z|)^{d-1} $, 
	\end{itemize}
	for all $z\in \mathbb{D} $ and $d\in \mathbb{N}. $ 
\end{lem}
\begin{proof} The first claim follows easily from the Schwarz-Pick lemma. In fact, by Schwarz-Pick we have 
		\[\frac{1+|\varphi(0)|}{1-|\varphi(0)|}(1-|\varphi(z)|^2) \geq \frac{(1-|\varphi(0)|^2)(1-|\varphi(z)|^2)}{|1-\overline{\varphi(0)}\varphi(z)|^2} =  1-\Big| \frac{\varphi(0)-\varphi(z)}{1-\overline{\varphi(0)}\varphi(z)} \Big|^2 \geq 1-|z|^2. \]
The second claim is from \cite[Theorem 1.1]{Dai2007}.
	\end{proof}

One final piece of information that we will use is related to the problem of similarity to a contraction. 
As already mentioned in the introduction, a polynomially bounded operator is not necessarily similar to a contraction. 
Nonetheless, a certain strengthening of the polynomial boundedness condition characterizes operators which are similar to contractions. 
We say that an operator  $T$ on a Hilbert space $H $  is completely polynomially bounded if there exists a constant $C $ such that  whenever $P = (p_{ij})_{i,j=1}^n $ is an $n\times n $ array of polynomials, the operator $P(T):=(p_{ij}(T))_{i,j=1}^n $ which acts on $ H \otimes \mathbb{C}^n  $ satisfies the inequality 
\[ \Vert P(T) \Vert \leq C \sup_{|z|<1} \Vert (p_{ij}(z))_{i,j=1}^n\Vert =:\Vert (p_{ij})_{i,j=1}^n\Vert. \]

In \cite{PAULSEN19841}, Paulsen gave the following characterization of operators which are similar to a contraction. 

\begin{thm}\label{thm:Paulsen} An operator on a Hilbert space $ H $ is similar to a contraction if and only if it is completely polynomially bounded.
	
\end{thm}

\section{Proof of main results} \label{sec:proofs}

In the sequel when $p $ is a holomorphic function in an open neighbourhood of the unit disc and $T $ is an operator such that $\sigma(T) \subseteq \overline{\mathbb{D}} $ by  $p(T) $ we will denote the Dunford-Riesz functional calculus of $T $. 
Observe that by approximating $p $ uniformly by polynomials in a disc of radius larger than $1 $,   equation \eqref{eq:FHpol} continues to hold for such functions $p $.  

\subsection{Proof of Theorems \ref{thm:FH_Kreiss} and \ref{thm:FHPb}}

The first part of this section is dedicated to the proof of Theorem \ref{thm:FH_Kreiss}. 

\begin{proof}[Proof of Theorem \ref{thm:FH_Kreiss}]

	We first rewrite the Kreiss condition in a more convenient form.
	For $\mu\in\mathbb D$ let
	\[
	R_\mu(z)=\frac{1}{1-\mu z}, z\in \mathbb{D}, \]
so that $ R_\mu(\Gamma_f)=(1-\mu \Gamma_f)^{-1}.$ 
By a straightforward calculation we also see that $R_\mu(S) $ is the multiplication operator $ M_{R_{\mu}} $. 
	
Since $R_\mu$ is analytic on the disc of radius $1/|\mu|$ and the spectrum of $\Gamma_f $ is contained in $\overline{\mathbb{D}} $, we obtain 

\begin{equation*}
	R_\mu(\Gamma_f) = \begin{bmatrix} R_\mu(S^*) & JH_{R'_{\mu}(S^*) f} \\ 0 & R_\mu(S) \end{bmatrix} = \begin{bmatrix} M_{R_{\overline{\mu}}}^* & JH_{M_{R'_{\overline{\mu}}}^* f}  \\ 0 & M_{R_\mu} \end{bmatrix}.
\end{equation*}

\noindent Therefore, 

\begin{align*}
	(1-|\mu|) \Vert R_\mu(\Gamma_f) \Vert &= (1-|\mu|) \Big\Vert  \begin{bmatrix} M_{R_{\overline{\mu}}}^* & JH_{M_{R'_{\overline{\mu}}}^* f}  \\ 0 & M_{R_\mu} \end{bmatrix} \Big\Vert \\ 
					    & = (1-|\mu|) \Big\Vert \begin{bmatrix} 0 & JH_{M_{R'_{\overline{\mu}}}^* f}  \\ 0 & 0 \end{bmatrix} + \begin{bmatrix} M_{R_{\overline{\mu}}}^* &   0 \\ 0 & M_{R_\mu} \end{bmatrix} \Big\Vert \\
					    & \geq (1-|\mu|) \Vert JH_{M_{R'_{\overline{\mu}}}^* f}  \Vert - (1-|\mu|) \max\{ \Vert R_{\overline{\mu}}\Vert_\infty, \Vert R_\mu\Vert_\infty \} \\
					    & = (1-|\mu|) \Vert JH_{M_{R'_{\overline{\mu}}}^* f}  \Vert - 1.
\end{align*}

Consequently, under the hypothesis that $ \Gamma_f $ satisfies the Kreiss condition we obtain the inequality
\begin{equation}\label{eq:Nehari_res} \sup_{|\mu|<1} (1-|\mu|) \Vert JH_{M_{R'_{\overline{\mu}}}^* f} \Vert < + \infty. \end{equation}
An application of Nehari's theorem then gives 
\begin{align*}
	\Vert JH_{M_{R'_{\overline{\mu}}}^* f}\Vert & = \Vert H_{M_{R'_{\overline{\mu}}}^* f}\Vert \\
						    & = \sup_{u\in H^2, \,\, \Vert u \Vert_1 \leq 1 } | \langle M_{R'_{\overline{\mu}}}^*f, u  \rangle| \\
						    & = \sup_{u\in H^2, \Vert u\Vert_1 \leq 1 } |\langle f, R'_{\overline{\mu}} u  \rangle| \\
						    & = \sup_{u\in H^2, \Vert u\Vert_1 \leq 1} \Big| \int_\mathbb{T} \frac{\mu f(z) \overline{u(z)}}{(1-\mu \overline{z})^2}\frac{|dz|}{2\pi} \Big|.
\end{align*}
In particular, since the functions $ u(z) = \frac{1-|\mu|^2}{(1-\overline{\mu}z)^2} $ have unit norm in $H^1 $ we conclude by \eqref{eq:Nehari_res} that $f $ must satisfy the inequality 
\[ \sup_{|\mu|<1} (1-|\mu|^2)^2 \Big| \int_{\mathbb{T}} \frac{\mu f(z)}{(1-\mu\overline{z})^4} \frac{|dz|}{2\pi}   \Big| =:C < +\infty. \]
This allows us to conclude the proof, because, since $f$ belongs to $H^2 $, so does the function $g(\mu) = \mu^3 f(\mu) $ and so it can be represented by the Cauchy integral of its boundary values

\[ g(\mu)= \int_{\mathbb{T}} \frac{z^3f(z)}{1-\mu \overline{z}}\frac{|dz|}{2\pi}.  \]
Which, after differentiating thrice in $\mu $ yields 
\begin{equation*}
	g'''(\mu)= 6 \int_\mathbb{T} \frac{f(z)}{(1-\mu \overline{z})^4} \frac{|dz|}{2\pi}.
\end{equation*}

In conjunction with the previous inequality we can conclude that $(1-|\mu|^2)^2|g'''(\mu)| \leq C  $ for all $\mu \in \mathbb{D} $, which is equivalent to the fact that $g'(z)= z^3f'(z) + 3z^2 f(z) $ belongs to the Bloch space \cite[Theorem 5.1.5]{Zhu1990}. Furthermore $\mathcal{B} $  is invariant under multiplication by $ z $ and $f\in BMOA $ which is a subspace of $\mathcal{B} $  which implies that $z \mapsto z^3f'(z) $ belongs to $\mathcal{B} $. Finally, the Bloch space is also invariant under the backwards shift operator $f \mapsto ( z\mapsto (f(z)-f(0))/z) $ which allows us to conclude that $f'\in \mathcal{B}. $

\end{proof}

Let us now turn our attention to the generalized Foguel-Hankel operators. 
The next lemma establishes the basic algebraic property of generalized Foguel-Hankel operators 

\begin{lem}\label{lem:FHAlg}

Let $f\in BMOA $ and $\varphi :\mathbb{D} \to \mathbb{D} $. If $p $ is holomorphic in an open neighbourhood of $\overline{\mathbb{D}}$ we have that
	\begin{equation}  
	p(\Gamma_{f,\varphi}) = \begin{bmatrix} M^*_{\widetilde{p\circ \varphi}} & JH_{M^*_{\widetilde{p'\circ \varphi}}f} \\ 0 & M_{p \circ \varphi} \end{bmatrix} = \Gamma_{M^*_{\widetilde{p'\circ \varphi}  }f, p\circ \varphi}.
\end{equation}
	
\end{lem}

\begin{proof}
	The proof is essentially based on the one for the classical Foguel-Hankel operators. 
	Notice that for $u,h\in H^2 $ 
	\begin{align*}
		\langle M^*_{\widetilde{\varphi}}JH_f u, h  \rangle  = \langle JH_fu	, \widetilde{\varphi} h   \rangle 
								 = \langle f , \widetilde{u \varphi} h \rangle 
								= \langle JH_fM_\varphi u, h \rangle.
	\end{align*}
Hence, 
\[ M_{\widetilde{\varphi}}^* J H_f = J H_f M_\varphi = JH_{M^*_{\widetilde{\varphi}}f}. \,\, \]
Since $\Gamma_{f,\varphi}$ is upper triangular, a computation immediately gives
	\[
	\Gamma_{f,\varphi}^n
	=
	\begin{bmatrix}
		M_{\widetilde{\varphi}}^{*n} &
		\displaystyle\sum_{k=0}^{n-1}
		M_{\widetilde{\varphi}}^{*k}\,JH_f\,M_\varphi^{\,n-1-k}\\[2mm]
		0 & M_\varphi^n
	\end{bmatrix} = 
	\begin{bmatrix}
		M^{*n}_{\widetilde{\varphi}} & n JH_f M_\varphi ^{n-1}\\
		0 & M_\varphi^n
	\end{bmatrix}.
	\]
The general case of a polynomial follows by linearity. Furthermore, for any polynomial $p $ we have the trivial norm estimate 
\[ \Vert \Gamma_{f,\varphi }\Vert \leq \Vert \varphi\Vert_\infty + \Vert f\Vert_*,  \]
which, since $p(\Gamma_{f,\varphi}) $ is another generalized Foguel-Hankel operator, gives that 
\[ \Vert p(\Gamma_{f,\varphi })\Vert \leq \sup_{z\in \varphi(\mathbb{D})}|p(z)| + \Vert M^*_{\widetilde{p'\circ \varphi}}f\Vert_* \leq  \sup_{z\in \varphi(\mathbb{D})}|p(z)| + \sup_{z \in \varphi(\mathbb{D})}|p'(z)| \Vert f\Vert_*. \]

Therefore if $p $ is a holomorphic function in an open neighbourhood of $ \overline{\mathbb{D}}  $ by approximating $p $ uniformly by polynomials in slightly larger disc we get that the conclusion of the lemma holds for $p $ as well.
\end{proof}

Before proceeding with the proof of Theorem \ref{thm:FHPb} we will need the following lemma. 

\begin{lem} \label{lem:Carphi}
	Let $ \varphi : \mathbb{D} \to \mathbb{D} $ be a holomorphic self-map of the unit disc and $r\in \mathbb{N}. $ There exists a constant $ C>0 $ depending on $\varphi  $ and $r $ such that for all arcs $I \subseteq \mathbb{T} $ and all $n \in \mathbb{N} $ 
	\begin{equation*}
		\frac{1}{|I|}\int_{S(I)}|\varphi(z)|^n(1-|\varphi(z)|)^rd m(z) \leq \frac{C}{n^{r+1}}.
	\end{equation*}
\end{lem}

\begin{proof}
Throughout the proof we will denote by $C $ a positive constant which depends only on $r $ and $\varphi $ and might change value from line to line. 	Consider the function $f(x)=x^n(1-x)^r, 0\leq x\leq 1 $ which has global maximum at $x_n^* =\frac{n}{r+n} $ and it is increasing in $[0,x_n^*] $. 
Furthermore, $f(x^*_n) \leq r^r/n^r $. 
	It follows from Lemma \ref{lem:selfmaps} (i) that there exists $k>0 $, depending only on $\varphi(0) $,  such that 
	\[ |\varphi(z)| \leq \frac{1+k|z|}{1+k}, |z|<1. \]
Consider the case that $ n > r/k$, so that $ x^*_n > \frac{1}{1+k} $. Let then $\rho_n \in (0,1) $ be the critical radius defined by solving the equation 
\[ \frac{1+k\rho_n}{1+k} = x^*_n. \]
By definition, $|\varphi(z)|\leq x^*_n $ whenever $ |z| \leq \rho_n. $ We therefore have that for an arc $ I \subseteq \mathbb{T} $ 

           \begin{align*} 
		   \int_{S(I)\cap \{|z| \leq \rho_n \}} f(|\varphi(z)|)dm(z) & \leq \int_{S(I)\cap \{ |z| \leq \rho_n \}}f\Big(\frac{1+k|z|}{1+k}\Big) d m(z) \\  
						     & \leq \frac{|I|}{\pi} \int_0^1 f\Big( \frac{1+k\rho}{1+k} \Big) d\rho \\
 & \leq \frac{|I| (k+1)}{\pi k} \int_0^1 x^n(1-x)^r  dx \\
 & = \frac{|I|(k+1)}{\pi k} B(n+1,r+1).
       \end{align*}
      In the last line $B $ is the Beta function. An application of Stirling's formula shows that there exists $C>0 $  such that $ B(n+1,r+1) \leq C/n^{r+1} $.  
The remaining part of the integral can be handled as follows 

     \begin{align*}
	     \int_{S(I) \cap \{ \rho_n \leq |z| <1 \}} f(|\varphi(z)|)d m(z) & \leq f(x_n^*) m(S(I)\cap \{ \rho_n \leq |z| < 1 \}) \\ 
							  & \leq f(x_n^*) \frac{|I|}{\pi}\Big( 1 - \frac{(1+k)x_n^*-1}{k} \Big) \\
							  & = \frac{|I|}{\pi}(k+1)f(x^*_n)(1-x^*_n) \\
							  & \leq C \frac{|I|}{n^{r+1}}.
     \end{align*}

Suppose finally that  $ x^*_n \leq \frac{1}{1+k} $ or equivalently $ n \leq r/k $.
Then, using the trivial estimate $ f\leq 1 $,  for any arc $ I\subseteq \mathbb{T} $ we have
\[ \int_{S(I)}f(|\varphi(z)|) dm(z) \leq (r/k)^{r+1} \frac{|I|}{n^{r+1}}.   \]
This concludes the proof of the lemma.
\end{proof}

\begin{lem}\label{lem:pb_equiv} Let $f\in BMOA $ and $\varphi : \mathbb{D}\to \mathbb{D} $ a holomorphic self-map of the unit disc. Then the operator $\Gamma_{f,\varphi} $ is power bounded if and only if
	\begin{equation}\label{eq:Nehari} \sup_{n \geq 1 } n \Vert JH_fM_{\varphi}^{n-1} \Vert = \sup_{n\geq 1} n \sup_{u\in H^2, \Vert u\Vert_1 \leq 1 } |\langle f, \widetilde{\varphi}\,^{n-1}u \rangle| < + \infty. 
	\end{equation}

\end{lem}

\begin{proof}
	Suppose that $f $ and $\varphi $ as in the statement. By virtue of Lemma \ref{lem:FHAlg}, showing that $\Gamma_{f,\varphi} $ is power bounded is equivalent to showing that   $  n J H_{M_{\widetilde{\varphi}}^{*n-1}f}  $ is uniformly bounded in norm. That is because the operators $M_{\widetilde{\varphi}}^*, M_\varphi $ are contractions.   
		
		Furthermore,  by Nehari's theorem 
		\[   \Vert JH_fM_{\varphi}^{n-1} \Vert = \Vert  J H_{M_{\widetilde{\varphi}}^{*(n-1)}f} \Vert = \sup_{u\in H^2, \Vert u\Vert_1 \leq 1 }| \langle M^{*(n-1)}_{\widetilde{\varphi}} f, u \rangle|  = \sup_{u\in H^2, \Vert u\Vert_1\leq 1}| \langle f, \widetilde{\varphi}\,^{n-1}u \rangle|.\]
Combined with the previous observation this equality concludes the proof of the Lemma.
\end{proof}
Before proving Theorem 2, we establish the following
geometric lemma, which will be used in its proof.

\begin{lem}\label{lem:carleson-localization}
	Let $I\subset\mathbb T$ be an arc. There exist absolute constants $C>0$ and $c>1$ such that
	\[
	\int_{S(I)}
	(1-|z|)^{-1}\mathbbm{1}_{D(z)}(w)\,dm(z)
	\leq
	C(1-|w|)\mathbbm{1}_{S(cI)}(w),
	\qquad w\in\mathbb D.
	\]
\end{lem}

\begin{proof} 
	Let $w\in\mathbb D$, we define $E(w):=\{z\in S(I): w\in D(z)\}.$ We first observe that
	\[
	\int_{S(I)}
	(1-|z|)^{-1}\mathbbm{1}_{D(z)}(w)\,dm(z)
	=
	\int_{E(w)}\frac{dm(z)}{1-|z|}.
	\]
	
	If $E(w)$ is empty there is nothing to prove. Suppose therefore that
	 there exists $z\in S(I)$ such that
	$w\in D(z)$. We then have that 
	\[
		2|z-w| < 1-|z| \leq 1-|w| +|z-w| .
	\]
	Rearranging the inequality we get that $E(w) \subseteq B(w,1-|w|):=\{ z\in \mathbb{D}: |w-z| < 1-|w| \} $. Moreover 
	\[ 1-|w| \leq 1-|z| + |z-w| < \frac{3}{2}(1-|z|). \]
		Consequently,
	\[
	\begin{aligned}
		\int_{E(w)}\frac{dm(z)}{1-|z|}
		&\leq
		\frac{3}{2}\frac{1}{1-|w|}\,m(E(w))\\
		&\leq
		\frac{3}{2}\frac{1}{1-|w|}
		\,m\bigl(B(w,1-|w|)\bigr)\\
		&\leq
		\frac{3\pi}{2}(1-|w|).
	\end{aligned}
	\]
	
	It remains to prove that the function at the left hand side of the inequality is localized in an enlarged Carleson box $S(c I) $. 
	For $z,w $ as before we have
	\begin{equation}\label{eq:elementary1}
	1-|w| \leq \frac{3}{2}(1-|z|) \leq \frac{3}{2}|I|.
	\end{equation}
	We now show that $w/|w|$ belongs to a fixed dilation
	of $I$. If $|I|\geq 1/2$, then of course  $16I=\mathbb T$, and hence $S(16I)=\mathbb D$. Thus  the desired conclusion is immediate.
	Now assume that $|I| \leq 1/2 $. 
	Since $z\in S(I)$,	
	\[
	|z|\geq 1-|I|>\frac{1}{2}.
	\]
	Furthermore,
	\[
	\begin{aligned}
		|w|\geq |z|-|w-z|>
		|z|-\frac{1-|z|}{2}=
		\frac{3|z|-1}{2}>\frac{1}{4}.
	\end{aligned}
	\]
Then
	\begin{equation*}
		\left|
		\frac{w}{|w|}
		-
		\frac{z}{|z|}
		\right|
		=
		\frac{||z|(w-z)+z(|z|-|w|)|}{|w||z|} 
		\leq
		8 (|w-z|+\bigl||w|-|z|\bigr|)
		\leq
		16|w-z|
	\end{equation*}
	Therefore, if $\zeta_I $ is the center of the arc $I $, using $C $ to denote a generic absolute constant, we have  
	\begin{equation}\label{eq:elementary2}
		\left| \frac{w}{|w|}-\zeta_I\right| \leq  \left| \frac{w}{|w|} -\frac{z}{|z|} \right| + \left| \frac{z}{|z|} - \zeta_I \right| \leq
	C|w-z| + C (1-|z|) \leq C(1-|z|) \leq C|I|.
	\end{equation}
	
	Since the Euclidean distance and the arc-length distance on the circle are locally comparable we conclude from inequalities \eqref{eq:elementary1} and \eqref{eq:elementary2} that there exists an constant $c>0 $ such that $w/|w| \in c I $. 
\end{proof}
We will now give the proof of Theorem \ref{thm:FHPb}. We will adopt the convention that the constant $C $ denotes a constant which may change from line to line, but it is independent of the natural number $n. $ 

	\begin{proof}[Proof of Theorem \ref{thm:FHPb}]

		Let $\psi =\widetilde{\varphi} $. Without loss of generality we may assume that $f(0)=f'(0)=0$. Therefore
		Lemma~\ref{lem:PL} reduces the estimate to an integral involving
		$f'' (\psi^n u)''$. Hence by expanding $(\psi^n u)'' $ we get
		\[
		\begin{aligned}
			(\psi^n u)''
			=
			n(n-1)\psi^{n-2}(\psi')^2u
			+n\psi^{n-1}\psi''u
			+2n\psi^{n-1}\psi'u'
			+\psi^n u''.
		\end{aligned}
		\]
		Using Lemma~\ref{lem:PL}, the estimate
		\begin{equation}\label{eq:kernelestimate}
			\log\frac1{|z|}(1+|z|^2)-(1-|z|^2)
			\leq
			\begin{cases}
				C(1-|z|)^3, & \frac12<|z|<1,\\
				C\log\frac1{|z|}, & 0<|z|<\frac12,
			\end{cases}
		\end{equation}
		and absorbing numerical constants into $C$, we obtain
		
		\begin{align*}
			n\Vert JH_{M^{*n}_\psi f}\Vert & \leq C \sup_{\Vert u\Vert_{1}\leq 1 } \Big(  n^3 \int_{\mathbb{D}} |\psi(z)|^{n-2}|\psi'(z)|^2|f''(z)| |u(z)|(1-|z|)^3 d m(z) \\
						       & \quad + n^2 \int_\mathbb{D} |\psi(z)|^{n-1} |\psi''(z)|  |u(z)||f''(z)|(1-|z|)^3 d m(z) \\ 
			&  \quad\quad  + n^2 \int_{\mathbb{D}} |\psi(z)|^{n-1}|\psi'(z)||u'(z)||f''(z)|(1-|z|)^3 d m(z)\\ 
			& \quad \quad \quad  + n \int_{\mathbb{D}}|\psi(z)|^n|u''(z)||f''(z)|(1-|z|)^3 d m(z) \Big) \\
			& \leq C \Vert f'\Vert_{\mathcal{B}} \sup_{\Vert u\Vert_{1} \leq 1 } \Big( n^3  \int_{\mathbb{D}}   |\psi(z)|^{n-2}|\psi'(z)|^2 |u(z)| (1-|z|)^2 d m(z) \\ 
			& \quad + n^2 \int_\mathbb{D} |\psi(z)|^{n-1} |\psi''(z)| |u(z)|(1-|z|)^2 d m(z) \\ 
			&  \quad \quad + n^2 \int_{\mathbb{D}} |\psi(z)|^{n-1}|\psi'(z)||u'(z)|(1-|z|)^2 d m(z)\\ 
			 &  \quad \quad \quad + n \int_{\mathbb{D}}|\psi(z)|^n|u''(z)|(1-|z|)^2 d m(z) \Big). 
		\end{align*} 
		
Let us call $I_k, k=1,2,3,4$ the four integrals appearing in the inequality above. 
We start by examining $I_1. $ Let $u\in H^1 $ such that $\Vert u\Vert_1 \leq 1$, and apply the Schwarz-Pick inequality to obtain
\begin{equation*}
	I_1 = \int_{\mathbb{D}}|\psi(z)|^{n-2}|\psi'(z)|^2|u(z)|(1-|z|)^2 d m(z) \leq C \int_{\mathbb{D}}|\psi(z)|^n (1-|\psi(z)|)^2 |u(z)| d m(z).
\end{equation*}
Then, an application of Lemma \ref{lem:Carphi} and Carleson's embedding theorem  gives the desired inequality, i.e. $I_1 \leq Cn^{-3} $. 
For $I_2$, Lemma~\ref{lem:selfmaps} $(ii)$ with $d=2$ gives
\[
|\psi''(z)|(1-|z|)^2\leq C(1-|\psi(z)|).
\]
Hence
\[
I_2
\leq
C\int_{\mathbb D}|\psi(z)|^{n-1}(1-|\psi(z)|)|u(z)|\,dm(z),
\]
and the same argument, using Lemma~\ref{lem:Carphi} with $r=1$, yields
\[
I_2\leq \frac{C}{n^2}.
\]

We now estimate $I_3$ by applying Theorem~\ref{thm:Lue} with $d=1$ to the
measure
\[
d\mu_n(z)=|\psi(z)|^{n-1}|\psi'(z)|(1-|z|)^2\,dm(z).
\]
Using again Schwarz-Pick and the comparability of $1-|w|$ and $1-|z|$ for
$w\in D(z)$, we get
\[
\mu_n(D(z))
\leq
C(1-|z|)
\int_{D(z)}|\psi(w)|^{n-1}(1-|\psi(w)|)\,dm(w).
\]
Consequently,
\[
\begin{aligned}
	\frac{\mu_n(D(z))^2}{(1-|z|)^5}
	&\leq
	C(1-|z|)^{-3}
	\left(\int_{D(z)}|\psi(w)|^{n-1}(1-|\psi(w)|)\,dm(w)\right)^2 \\
	&\leq
	C(1-|z|)^{-1}
	\int_{D(z)}|\psi(w)|^{2n-2}(1-|\psi(w)|)^2\,dm(w).
\end{aligned}
\]

Therefore, using Lemma \ref{lem:carleson-localization}
\[
\begin{aligned}
	\frac1{|I|}\int_{S(I)}
	\frac{\mu_n(D(z))^2}{(1-|z|)^5}\,dm(z)
	&\leq
	\frac{C}{|I|}\int_{S(cI)}
	(1-|w|)|\psi(w)|^{2n-2}(1-|\psi(w)|)^2\,dm(w) \\
	&\leq
	\frac{C}{|I|}\int_{S(cI)}
	|\psi(w)|^{2n-2}(1-|\psi(w)|)^3\,dm(w) \\
	&\leq \frac{C}{n^4},
\end{aligned}
\]
where in the second inequality we used Lemma~\ref{lem:selfmaps}, and in the last
one Lemma~\ref{lem:Carphi} with $r=3$. By Theorem~\ref{thm:Lue}, this implies
\[
I_3\leq \frac{C}{n^2}.
\]

Finally, $I_4$ is estimated in the same way, now applying Theorem~\ref{thm:Lue}
with $d=2$ to the measure
\[
d\nu_n(z)=|\psi(z)|^n(1-|z|)^2\,dm(z).
\]
Indeed,
\[
\nu_n(D(z))
\leq C(1-|z|)^2\int_{D(z)}|\psi(w)|^n\,dm(w),
\]
and hence
\[
\frac1{|I|}\int_{S(I)}
\frac{\nu_n(D(z))^2}{(1-|z|)^7}\,dm(z)
\leq
\frac{C}{|I|}\int_{S(cI)}(1-|w|)|\psi(w)|^{2n}\,dm(w).
\]
Using again Lemma \ref{lem:carleson-localization}, Lemma~\ref{lem:selfmaps} and then Lemma~\ref{lem:Carphi} with
$r=1$, we get
\[
\frac1{|I|}\int_{S(I)}
\frac{\nu_n(D(z))^2}{(1-|z|)^7}\,dm(z)
\leq \frac{C}{n^2}.
\]
Thus Theorem~\ref{thm:Lue} gives
\[
I_4\leq \frac{C}{n}.
\]

Combining the estimates for $I_1,I_2,I_3,I_4$, we obtain that $\Gamma_{f,\varphi}$ is power bounded.	\end{proof}

	\subsection{Power boundedness for a fixed symbol $\varphi $. } \label{sec:fixed_symbol}

Theorem \ref{thm:FHPb} provides a sufficient condition for the generalized Foguel-Hankel operator to be power bounded. 
In contrast to the case of classical Foguel-Hankel operators this condition is no longer necessary in general. For example this can be seen by examining strictly contractive symbols $\varphi $. Consider a holomorphic function $\varphi : \mathbb{D} \to \mathbb{D} $ such  that $ \Vert \varphi\Vert_{\infty}<1 $, $f\in BMOA $ and $u\in H^{2}. $ 
Then we have 
\[ n|\langle f, \tilde{\varphi}^{n-1} u \rangle| \leq n \Vert f\Vert_* \Vert \varphi\Vert_\infty^{n-1}\Vert u\Vert_1 \]
and hence, by Lemma \ref{lem:pb_equiv}, $\Gamma_{f,\varphi} $ is power bounded. 
Concretely, taking  

\[ f(z) = \log \frac{1}{1-z} \]
and $\varphi(z)=z/2 $ we have that $ f\in BMOA $ but $f'\not\in \mathcal{B} $ because 
\[ (1-r^2)|f''(r)| = \frac{1+r}{1-r} \to \infty. \]
Therefore the operator $\Gamma_{f,\varphi} $ is power bounded even though $f'\not\in \mathcal{B}. $
It is also instructive to consider the case that $\varphi(z) = z^m $ for some $m\in \mathbb{N} $, in other words $ M_{\varphi}  = S^m. $  

\begin{prop}
	Let $ f\in BMOA $ and $ e_m(z) = z^m $ for some $m\in \mathbb{N} $. Then the operator $ \Gamma_{f,e_m}$ is power bounded if and only if $f'\in \mathcal{B} $.  
\end{prop}

\begin{proof}
The sufficiency in the condition is contained in Theorem \ref{thm:FHPb}.  By Lemma \ref{lem:pb_equiv} $\Gamma_{f,e_m} $ is power bounded if and only if 
	\[ \sup_{n\geq 1} n \Vert JH_f S^{m(n-1)}\Vert < + \infty. \]
Let   
\[ a_N := \Vert JH_f S^N \Vert. \]
Therefore, there exists a constant $C>0 $ such that $ n a_{m(n-1)} \leq C $ for all $ n\geq 1 $.  Notice that the sequence $(a_N) $ is non-increasing, since $ S^{N+1}H^2 \subseteq S^N H^2 $. Suppose now that $ N \geq 2 $ and it satisfies $ m(n-1) \leq N < mn $, then by monotonicity 
\[ (N+1) a_N \leq (nm+1) a_{m(n-1)} \leq (m+1) C. \]
Hence, using again the necessity part of Lemma \ref{lem:pb_equiv}, we arrive at the conclusion that $ \Gamma_{f,e_1} = \Gamma_f $ is power bounded, which via Peller's condition \cite{Peller1982} implies that $f' \in \mathcal{B}. $ 

\end{proof} 

The next natural class of symbols $\varphi $ that one would like to consider besides strict contractions and monomials would be Blaschke products and singular inner functions. We expect that in the case of finite Blaschke products the condition of Peller should continue to be necessary for the power boundedness of the operator $\Gamma_{f,\varphi} $ although we do not have a proof of this claim. The case of singular inner functions is possibly even more intricate and we do not have a plausible conjecture for the characterization in this case. 

\section{Hilbert-Foguel operators} \label{sec:Hilbert}

In what follows, the operator $JH_f $ will be the Hilbert matrix $\mathcal{H}$, which is represented, with respect to the standard orthonormal basis $H^2 $,  by the matrix 
\begin{equation*}
	\mathcal{H} = \Big[ \frac{1}{n+m+1} \Big]_{n,m=0}^\infty.
\end{equation*}

With a slight abuse of notation we will denote by the same symbol the corresponding operator on $H^2$.
Therefore, given $u,v\in H^2 $ we have that 

\begin{align*}
	\langle \mathcal{H}u, v \rangle & = \sum_{n=0}^{\infty} \sum_{m=0}^{\infty} \frac{\widehat{u}(m)\overline{\widehat{v}(n)} }{n+m+1}
\\ 
					& = \int_0^1 \sum_{n=0}^{\infty} \sum_{m=0}^{\infty} \widehat{u}(m)\overline{\widehat{v}(n)}t ^{n+m} dt \\
					& = \int_0^1 u(t)\overline{v(t)}dt.
\end{align*}

\begin{proof}[Proof of Theorem \ref{thm:Hilbert}]

	We will first prove that $(iii)\implies (i)$. 
	This is a direct application of Theorem \ref{thm:Paulsen}. 
	Observe that the simple argument in \cite[p. 153]{Davidson1997} does not apply in this case. 
	Let $(p_{ij})_{i,j=1}^n $ be an $n\times n $ matrix of analytic polynomials $n=1,2,\dots $ and let $(x_i)_{i=1}^n, (y_j)_{j=1}^n \in H^2 \oplus H^2 $. 
	What we need to show is that there exists a constant $C>0 $ such that 
	\begin{equation}
		\Big| \sum_{i,j=1}^{n} \langle p_{ij}(h_\varphi) x_i, y_j \rangle \Big| \leq C \Vert (p_{ij})_{i,j=1}^n\Vert \Big( \sum_{i=1}^{n}  \Vert x_i\Vert^2 \Big)^{\frac{1}{2}}\Big( \sum_{j=1}^{n} \Vert y_j\Vert^2 \Big)^{\frac{1}{2}}.
	\end{equation}
	Let 
	\begin{equation*}
		x_i = 
		\begin{bmatrix}
			u_i \\ v_i
		\end{bmatrix}, 
		y_j = \begin{bmatrix}
			h_j \\ g_j
		\end{bmatrix}, \,\,\,  u_i,v_i,h_j,g_j \in H^2, i,j=1,\dots n,
	\end{equation*}
	such that 
\begin{equation*}
	 \sum_{i=1}^{n}  \Vert x_i\Vert^2 \leq 1, \,\,\, \sum_{j=1}^{n} \Vert y_j\Vert^2 \leq 1.
\end{equation*}	

\noindent We have that 
\begin{align*}
	\Big| \sum_{i,j=1}^{n} \langle p_{ij}(h_\varphi) x_i, y_j \rangle \Big| &  \leq \Big| 	\sum_{i,j=1}^{n} 
											  \langle
											 p_{ij} \Big( \begin{bmatrix}
												  M_{\widetilde{\varphi}}^* & 0 \\
												  0 & M_{ \varphi}
											  \end{bmatrix} \Big) x_i
											  ,  y_j \rangle \Big| + \\
											  & \,\,\,\,\, \Big| \sum_{i,j=1}^{n} \langle \mathcal{H} M_{p_{ij}'\circ \varphi} v_i , h_j  \rangle \Big| \\
											  & \leq \Vert (p_{ij})_{i,j=1}^n \Vert+ \Big| \sum_{i,j=1}^{n} \int_0^1 p'_{ij}(\varphi(t))v_i(t)\overline{h_j(t)}dt \Big| .
\end{align*}

In the above inequalities we have used that the diagonal operator appearing above is a contraction and therefore completely polynomially bounded with constant $1 $.

Now notice that by the hypothesis and the maximum principle there must exist $c<1 $ such that $ \sup_{0<t<1}|\varphi(t)|\leq c. $ 
Hence by Cauchy's formula for vector valued holomorphic functions we get 
\begin{equation*}
	\sup_{t\in (0,1)}\Vert (p_{ij}'\circ\varphi(t))_{i,j=1}^n \Vert \leq \sup_{|z|<c} \Vert (p_{ij}'(z))_{i,j=1}^n \Vert  \leq C \sup_{|z|<1} \Vert (p_{ij}(z))_{i,j=1}^n\Vert, 
\end{equation*}
for some constant $C $ independent of $(p_{ij})_{i,j=1}^n.$ 
With this information we can estimate the second sum as follows using the Fej\'er-Riesz inequality  \cite[Theorem 3.13]{Duren1970}
\begin{align*}
	\Big| \sum_{i,j=1}^{n} \int_0^1 p'_{ij}(\varphi(t))v_i(t)\overline{h_j(t)}dt \Big| & \leq C \Vert (p_{ij})_{i,j=1}^n\Vert  \int_0^1 \Big( \sum_{j=1}^{n} |v_i(t)|^2 \Big)^{\frac{1}{2}} \Big( \sum_{j=1}^{n} |h_j(t)|^2 \Big)^{\frac{1}{2}} dt   \\  
	& \leq C \Vert (p_{ij})_{i,j=1}^n\Vert \Big( \sum_{i=1}^{n} \int_0^1 |v_i(t)|^2 dt   \Big)^{\frac{1}{2}} \Big( \sum_{j=1}^{n} \int_0^1 |h_j(t)|^2 dt   \Big)^{\frac{1}{2}} \\
	 & \leq C \Vert (p_{ij})_{i,j=1}^n\Vert \Big( \sum_{i=1}^{n} \Vert v_i\Vert_2^2  \Big)^{\frac{1}{2}} \Big( \sum_{j=1}^{n} \Vert h_j\Vert_2^2 \Big)^{\frac{1}{2}} \\
	 & \leq C \Vert (p_{ij})_{i,j=1}^n\Vert.
\end{align*}

To conclude the proof, we show that $(ii) \implies (iii) $. Suppose that $(ii) $  holds. As it was the case in the proof of Theorem \ref{thm:FH_Kreiss} this implies that 
\begin{equation} \label{eq:test}
\sup_{|\mu|<1}(1-|\mu|) \Vert \mathcal{H}M_{R_{\overline{\mu}}'\circ \varphi }\Vert = \sup_{|\mu|<1} \sup_{\Vert u\Vert_2 \leq 1, \Vert v\Vert_2 \leq 1 }(1-|\mu|) \Big|  \int_0^1 \frac{\mu u(t) \overline{v(t)}}{(1-\mu \varphi(t))^2} dt \Big|  < + \infty. \end{equation}
Next, fix $\delta\in(0,1)$ and set $a=1-\delta$. For $|\mu|<1$ consider the analytic functions
\begin{equation}	u(z) = \frac{(1-a^2)^{\frac{1}{2}}}{1-az}, \,\,\, z\in \mathbb{D}, \, \, \, v(z)=\overline{\mu}\,
	\frac{(1-|\mu|)^{2}}{\bigl(1-\overline{\mu}\,\varphi(z)\bigr)^{2}} \frac{(1-a^2)^\frac{1}{2}}{1-az}.
	\end{equation}

	The fact that $\varphi $ is bounded in modulus by $1 $ implies that the functions $u, v $ have $H^2 $ norm less than $1 $ and so they are admissible for 
	\eqref{eq:test}, in other words there exists $C>0 $ such that for all $\mu \in \mathbb{D} $ and $\delta \in (0,1) $     
	
	\[
	\int_{0}^{1}
	\frac{|\mu|^{2}(1-|\mu|)^{3}}{|1-\overline{\mu}\,\varphi(\rho)|^{4}}\,
	\frac{1-a^{2}}{(1-a\rho)^{2}}\,d\rho
	\le C.
	\]
	Restricting to $\rho\in[1-\delta,1]$ and observing that on this interval
	\[
	1-a\rho=1-(1-\delta)\rho=(1-\rho)+\delta\rho\le 2\delta,
	\qquad
	1-a^{2}=\delta(2-\delta)\ge \delta,
	\]
	we get $\dfrac{1-a^{2}}{(1-a\rho)^{2}}\ge \dfrac{1}{4\delta},\, \rho\in[1-\delta,1]. $
	Therefore,
	\[
	\frac{(1-|\mu|)^{3}|\mu|^{2}}{4\delta}\int_{1-\delta}^{1}\frac{d\rho}{|1-\overline{\mu}\,\varphi(\rho)|^{4}}
	\le C, \,\, \text{for all} \,\, \delta \in (0,1), \mu \in \mathbb{D}.
	\]
		In particular, for $|\mu|\in[1/2,1)$ and $\delta\in(0,1/2]$ we obtain the uniform bound
	\begin{equation}\label{eq:delta_estimate_uniform}
		\sup_{\substack{|\mu|\in[1/2,1)\\ 0<\delta\le 1/2}}
		\frac{(1-|\mu|)^{3}}{\delta}\int_{1-\delta}^{1}\frac{d\rho}{|1-\overline{\mu}\,\varphi(\rho)|^{4}}
		<\infty.
	\end{equation}
	Assume, towards a contradiction, that
	\begin{equation}\label{eq:bad_limsup}
		\limsup_{\substack{t\to 1^-\\ t\in\mathbb R}}|\varphi(t)|=1.
	\end{equation}
	Then there exists $t_k\to 1$ such that $|\varphi(t_k)|\to 1$. Let
	\[
	\tau_k(z):=\frac{z+t_k}{1+t_k z},
	\qquad \tau_k(0)=t_k.
	\]
	The family $\{\varphi\circ\tau_k\}_{k\ge 1}$ is normal in $\mathbb{D}$; hence, after passing to a subsequence, there exists an analytic function $g$ on $\mathbb{D}$ such that
	\[
	\varphi\circ\tau_k \longrightarrow g
	\qquad\text{uniformly on compact subsets of } \mathbb{D},
	\qquad |g|\le 1.
	\]
	Moreover,
	\[
	|g(0)|
	=
	\lim_{k\to\infty}|\varphi(\tau_k(0))|
	=
	\lim_{k\to\infty}|\varphi(t_k)|
	=
	1,
	\]
	and thus the maximum modulus principle implies that $g$ is equal to a unimodular constant $e^{i\theta}, \theta \in \mathbb{R}. $ 
	Fix $r\in(1/2,1)$ and set $\mu:=r e^{i\theta}$, so $|\mu|=r$. Define
	\[
	\Phi_\mu(z):=\frac{1}{|1-\overline{\mu}\,\varphi(z)|^{4}}.
	\]
	Then $\Phi_\mu\circ\tau_k \to (1-r)^{-4}$ uniformly on compact subsets of $\mathbb{D}$, hence there exists $k(r) \in \mathbb{N} $, depending on $r $ such that for $k \geq k(r) $ 
	\[
	\Phi_\mu(\tau_k(s))\ge \frac{1}{2(1-r)^{4}}
	\qquad\text{for all } s\in\Bigl[0,\frac12\Bigr].
	\]
	Let $\delta_k:=1-t_k$, so that 	
	\begin{align*}
		\frac{(1-r)^{3}}{\delta_k}\int_{1-\delta_k}^{1}\Phi_\mu(\rho)\,d\rho
		&\ge
		\frac{(1-r)^{3}}{\delta_k}\int_{t_k}^{\tau_k(1/2)}\Phi_\mu(\rho)\,d\rho \notag\\
		&=
		\frac{(1-r)^{3}}{\delta_k}\int_{0}^{1/2}\Phi_\mu(\tau_k(s))\,\tau_k'(s)\,ds \notag\\
		&\ge
		\frac{(1-r)^{3}}{\delta_k}\cdot \frac{1}{2(1-r)^{4}}
		\int_{0}^{1/2}\tau_k'(s)\,ds \notag\\
		&=
		\frac{1}{2(1-r)}\cdot\frac{\tau_k(1/2)-t_k}{\delta_k}.
		\label{eq:lower_chain}
	\end{align*}
	A direct computation gives
	\[
	\tau_k\Bigl(\frac12\Bigr)-t_k
	=
	\frac{1-t_k^{2}}{2\left(1+\frac{t_k}{2}\right)}
	=
	\frac{(1-t_k)(1+t_k)}{2\left(1+\frac{t_k}{2}\right)}
	\ge
	\frac{1-t_k}{3}
	=
	\frac{\delta_k}{3},
	\]
	and therefore
	\begin{equation}\label{eq:lower_bound_final}
		\frac{(1-r)^{3}}{\delta_k}\int_{1-\delta_k}^{1}\frac{d\rho}{|1-\mu\,\varphi(\rho)|^{4}}
		\ge
		\frac{1}{6(1-r)}.
	\end{equation}
	For $k$ large we have $\delta_k\le 1/2$, so \eqref{eq:lower_bound_final} contradicts \eqref{eq:delta_estimate_uniform} by letting $r\to 1^-$, obtaining a contradiction. This shows that \eqref{eq:bad_limsup} is false.
\end{proof}

\section{Open problems} \label{sec:open_problems}

	It is quite evident that there are many questions that are natural in this context and to which we have not been able to give a satisfactory answer. 
	For example, in view of Theorem \ref{thm:FHPb} it is reasonable to conjecture that for every $f $ such that $f'\in BMOA $ and every holomorphic self-map of the unit disc $\varphi  $ the operator $\Gamma_{f,\varphi} $ is similar to a contraction. 
	Following \cite{Peller1982} it can be shown that this is the case if the operator 
\[ V_\varphi : H^\infty \to BMOA, \,\,\, V_\varphi f (z) = \int_0^z f'(\varphi(t)) dt  \]
is bounded, but we do not know if this holds for all $\varphi $.   
A more ambitious goal is to characterize all pairs $(f, \varphi) $ such that $\Gamma_{f,\varphi} $ is polynomially bounded, power bounded or satisfies the Kreiss condition.

In particular it would be interesting to know if there exists a  couple $f,\varphi $  such that $\Gamma_{f,\varphi} $ is Kreiss but not power bounded. 
Because of the upper triangular form of $ \Gamma_{f,\varphi} $, by rescaling $f $ if necessary, one would have an example of an operator with Kreiss constant arbitrarily close to $1 $ which is not power bounded. 
Such operators have been constructed in \cite{McCarthy1965, Chalmoukis2025}, but the quantitative comparison of the Kreiss condition and the power boundedness condition is not yet understood in the case when the Kreiss constant is close to $1. $ 

\section*{Acknowledgements}
The authors would like to thank the anonymous referee for a careful and thorough reading of the manuscript as well as suggesting numerous improvements and pointing out inaccuracies in the first version. In particular the discussion around Theorem \ref{thm:FHPb} in Section \ref{sec:fixed_symbol} was motivated by the reviewer's comments.

\bibliographystyle{amsplain}
\bibliography{biblio}

\providecommand{\bysame}{\leavevmode\hbox to3em{\hrulefill}\thinspace}
\providecommand{\MR}{\relax\ifhmode\unskip\space\fi MR }
\providecommand{\MRhref}[2]{%
  \href{http://www.ams.org/mathscinet-getitem?mr=#1}{#2}
}
\providecommand{\href}[2]{#2}
\begin{thebibliography}{10}

\bibitem{Aleksandrov1996}
A.~B. Aleksandrov and V.~V. Peller, \emph{Hankel operators and similarity to a contraction}, Int. Math. Res. Not. \textbf{1996} (1996), no.~6, 263--275 (English).

\bibitem{Bourgain1986}
J.~Bourgain, \emph{On the similarity problem for polynomially bounded operators on {Hilbert} space}, Isr. J. Math. \textbf{54} (1986), 227--241 (English).

\bibitem{Chalmoukis2025}
Nikolaos Chalmoukis, Georgios Tsikalas, and Dmitry Yakubovich, \emph{Operators with small {Kreiss} constants}, Preprint, {arXiv}:2512.10025 [math.{FA}] (2025), 2025, To appear in J. Anal. Math.

\bibitem{Dai2007}
Shaoyu Dai and Yifei Pan, \emph{Note on {S}chwarz-{P}ick estimates for bounded and positive real part analytic functions}, Proc. Am. Math. Soc. \textbf{136} (2008), no.~2, 635--640 (English).

\bibitem{Davidson1997}
Kenneth~R. Davidson, \emph{Polynomially bounded operators, a survey}, Operator algebras and applications. Proceedings of the NATO Advanced Study Institute and Aegean conference, Pythagorio, Samos, Greece, August 19--28, 1996, Dordrecht: Kluwer Academic Publishers, 1997, pp.~145--162 (English).

\bibitem{Duren1970}
P.~L. Duren, \emph{Theory of {{\({H}^ p\)}} spaces}, Pure Appl. Math., Academic Press, vol.~38, New York {and} London: Academic Press, 1970 (English).

\bibitem{Garnett2006}
John~B. Garnett, \emph{Bounded analytic functions}, revised 1st ed. ed., Grad. Texts Math., vol. 236, New York, NY: Springer, 2006 (English).

\bibitem{Luecking1991}
Daniel~H. Luecking, \emph{Embedding derivatives of {Hardy} spaces into {Lebesgue} spaces}, Proc. Lond. Math. Soc. (3) \textbf{63} (1991), no.~3, 595--619 (English).

\bibitem{McCarthy1965}
C.~A. McCarthy and J.~Schwartz, \emph{On the norm of a finite {Boolean} algebra of projections, and applications to theorems of {Kreiss} and {Morton}}, Commun. Pure Appl. Math. \textbf{18} (1965), 191--201 (English).

\bibitem{Neumann}
John~Von Neumann, \emph{Eine spektraltheorie für allgemeine operatoren eines unitären raumes. erhard schmidt zum 75. geburtstag in verehrung gewidmet}, Mathematische Nachrichten \textbf{4} (1950), no.~1-6, 258--281.

\bibitem{PAULSEN19841}
Vern~I Paulsen, \emph{Every completely polynomially bounded operator is similar to a contraction}, Journal of Functional Analysis \textbf{55} (1984), no.~1, 1--17.

\bibitem{Peller1982}
Vladimir~V. Peller, \emph{Estimates of functions of power bounded operators on {Hilbert} spaces}, J. Oper. Theory \textbf{7} (1982), 341--372 (English).

\bibitem{Peller1983}
\bysame, \emph{Estimates of functions of {H}ilbert space operators, similarity to a contraction and related function algebras}, Lect. Notes Math. \textbf{1573} (1984), 298--302 (English).

\bibitem{Pisier1996}
Gilles Pisier, \emph{A polynomially bounded operator on {Hilbert} space which is not similar to a contraction}, J. Am. Math. Soc. \textbf{10} (1997), no.~2, 351--369 (English).

\bibitem{Zhu1990}
Kehe Zhu, \emph{Operator theory in function spaces}, Pure Appl. Math., Marcel Dekker, vol. 139, New York etc.: Marcel Dekker, Inc., 1990 (English).

\end{thebibliography}

\end{document}